\documentclass[10pt]{article}

\usepackage[all]{xy}
\usepackage{amssymb,hyperref}
\usepackage{graphicx}
\usepackage{amsmath}
\usepackage{array}
\usepackage{indentfirst}

\newcommand{\F}{{\mathbb F}_2}

\newcommand{\Ac}{\mathcal{A}}

\newcommand{\gl}{\mathrm{GL}}

\newcommand{\T}{\mathbf T}

\def\L{\mathbf L}

\newcommand{\reg}{\widetilde{reg}}
\def\M{{\mathbf M}}

\newcommand{\U}{\mathcal U}
\newcommand{\Z}{\mathbb Z}
\newcommand\ra{\longrightarrow}

\newcommand\C{{\mathcal C}}

\newcommand\R{{\mathbb R}}
\newcommand\TT{{\mathbf T}}

\newcommand\Sq{\mathrm{Sq}}

\newcommand{\bqn}{\begin{eqnarray*}}
\newcommand{\eqn}{\end{eqnarray*}}

\author{Nguyen Dang Ho Hai, Lionel Schwartz}
\title{Realizing a complex of unstable modules}

\newtheorem{theorem}{Theorem}[section]
\newtheorem{proposition}[theorem]{Proposition}

\newtheorem{conjecture}[theorem]{Conjecture}

\newenvironment{remark}%
    {
	\begin{trivlist}
	\refstepcounter{theorem}%
    	\item[]{\textbf{Remark\ \thetheorem\ }}}%
    	{\end{trivlist}
    }

\begin{document}
\maketitle

\begin{abstract}
In a preceding article the authors and  Tran Ngoc Nam constructed a minimal injective resolution of the mod $2$
cohomology of a Thom spectrum.  A Segal conjecture type theorem for this spectrum was proved. In this paper one shows that the above mentioned resolutions can be realized topologically. In fact
there exists a family of cofibrations inducing short exact sequences in mod $2$ cohomology.  The  resolutions above are obtained 
 by splicing together these short exact sequences. Thus the injective resolutions are realizable in the best possible sense. In fact our construction
appears to be in some sense an injective closure of one of Takayasu. It strongly suggests that one can construct  geometrically (not only homotopically) certain dual Brown-Gitler spectra.
\end{abstract}

\tableofcontents


\section{Introduction}\label{introduction}

Below $H^*X$, resp. $H_*X$,  will denote the singular mod $2$ cohomology, 
resp. homology, of a space or a spectra  $X$. All spaces or spectra  will be supposed to be 
$2$-completed and to have finite mod $2$ cohomology in each degree. 
Let $\Ac$ be the mod $2$ Steenrod algebra and let $\mathcal U$ be the category of unstable $\Ac$-modules.

Let $L_n$ be the Steinberg summand in the mod $2$ cohomology $H^*B(\Z/2)^n$, as defined by Mitchell and Priddy \cite{MP83}. Let also $J(\ell)$ be the $\ell$-th Brown-Gitler module. These are both injective unstable $\Ac$-modules as well as their tensor product $L_k \otimes J(\ell)$ following Lannes and Zarati \cite{LZ86}.

The Brown-Gitler modules are the mod $2$ cohomology of spectra, even of spaces if $\ell$ is odd. We will denote these spectra by $\TT(\ell)$. If $\ell$ is odd this is a suspension spectrum.

Let $V_n=(\Z/2)^n$. Let us denote by  $L'_n$ the unstable module which is obtained as follows. Let $e_n$ be the Steinberg idempotent  in $\F[ \gl_n(\F)]$. Let 
$\reg_n\colon V_n \rightarrow O(2^n-1)$ be the regular reduced representation of  $V_n$ over the field $\R$.
 
Let  $\reg_n^{\oplus k}=\reg_n+\cdots+\reg_n$,  be the direct sum of $k$ copies of this representation. Let 
$BV_n^{\reg_n^{\oplus k}}$ be the Thom space associated to  $\reg_n^{\oplus k}$, {\it i.e.} the Thom space of the vector bundle $EV_n\times_{V_n}\R^{k(2^n-1)}\rightarrow BV_n$. The Steinberg idempotent determines a map of 
these Thom spaces after one suspension. Let us denote \cite{Tak99} by $\M(n)_k$  the Thom spectrum $e_n\cdot BV_n^{\reg_n^{\oplus k}}$. In particular if $k=1$ it  is the spectrum $\L(n) $ of  Mitchell and Priddy, and we will simply denote it $\L'(n)$ if
$k=2$. We will denote  by $L_n$ (as usual) and by $L'_n$ their mod $2$ cohomology.

Let us now recall:

\begin{theorem}[\cite{NST09b,NST09a}] \label{NST}
For all $n\ge 1$ there exists an acyclic complex of unstable injective 
 $\Ac$-modules
$$ 0 \rightarrow  L'_n\rightarrow
L_n\to L_{n-1}\otimes J(1)\to L_{n-2}\otimes J(3)\rightarrow \cdots\rightarrow
L(1)\otimes J(2^{n-1}-1)\rightarrow J(2^n-1)\rightarrow 0.
$$ 
This complex is a minimal injective resolution of $L'_n$.
\end{theorem}

Here is the main result of this article.

\begin{theorem}\label{rea}
The above complexes are realizable in the sense there exists a complex of spectra
$$
\TT(2^n-1) \ra \L(1) \wedge \TT (2^{n-1}-1)\ra \ldots \ra \L(i) \wedge \TT( 2^{n-i}-1) \ra \ldots \ra \L(n) \ra \L'(n)  
$$
so that when applying mod $2$ singular cohomolgy one gets the algebraic complexes. The composite of two successive maps is homotopically trivial, moreover all Toda brackets vanish.
\end{theorem}

The result follows from the following stronger result:

\def\R{\bf R}

\begin{theorem}\label{NS}
There are cofibrations of spectra, $1 \leq k \leq n$,
$$
\C(k-1) \ra \L(k)\wedge \TT(2^{n-k}-1)  \ra \C(k)
$$
so that 
\begin{itemize}
\item The induced sequences in mod $2$ cohomology are short exact;
\item $\C(n)=\L'(n)$;
\item $\C(0)=\TT(2^n-1)$.
\end{itemize}
\end{theorem}

Because of  \cite{NST09b,NST09a}, and further applications we are focused on the case of \ref{rea}. However this theorem is quite general. The proof of the general case proceeds exactly as in the particular case. Thus, we will only state:

\begin{theorem}\label{reag}
Let $M$ be the mod $2$ cohomology of a spectrum $X$. Let us assume that $M$ has an injective resolution of finite type in the category of unstable modules $\U$. By finite type one means that any term in the injective resolution is a finite direct sum of indecomposable modules. Then this resolution can be realized topologically by spectra, in the sense described above as soon as the algebraic map from $M$ to its injective hull in $\U$ can be realized.
\end{theorem}

The above theorems are in the category of spectra, however  the spectra $\L(i) \wedge \TT(2^{n-i}-1)$ are suspension spectra if $n>i$. In Theorem \ref{rea} only $\L(n)$ and $\L'(n)$ are not. They are only
direct summand of a suspension spectrum (after one suspension). It is thus reasonable to ask whether or not the complex of spectra can be realized as a complex of spaces (in an appropriate sense) after one suspension. The technology to do that is provided by the work of P. Goerss, J. Lannes and F. Morel \cite{GLM92}. However even if some preliminary computations are positive they are very technical in nature and we have not been able up to now
to decide this question. Thus we only propose:

\begin{conjecture} The preceding cofibrations in \ref{rea} can be realized as cofibrations of spaces, except for the last one that can be realized after one suspension and which is a retract of the suspension of a cofibration of spaces.
\end{conjecture}

In the general case it does not look plausible to hope for such a result, at least it is necessary to require something on the map realizing the embedding in the injective hull. Next we mentioned applications of \ref{rea}. It is clear that \ref{NS} leads to the construction of spectral sequences that we will consider either in a later version of this paper or in another one.


\section{Spectra satisfying the Brown-Gitler property and related spectra}

In this section we recollect some facts from \cite{HK00}. Most (but not all) of the material in this paper comes from earlier work, nicely collected here. Moreover we will need an original result of this paper. We shall say that a spectrum is {\it space-like} 
 if it is a wedge summand of a suspension 
spectrum. For such a spectrum $S$ the evaluation map  $ \Sigma^\infty \Omega^\infty S \ra S$ has a section and induces a monomorphism in cohomology. 

Let now $S$ be a spectrum whose mod $2$ singular cohomology is an injective unstable $\Ac$-module. Following  N. Kuhn we shall say that $S$ has the Brown-Gitler property if the natural map
$$
f \mapsto f^*,  \; \; \; \;   [ S, X]  \ra  {\rm Hom}_{\mathcal U}(H^*X,H^*S)
 $$
from homotopy classes from $S$ to $X$ to the set of maps of unstable modules  ${\rm Hom}_{\mathcal U}(H^*X,H^*S)$ is onto for all space-like spectrum $X$.

\begin{theorem}[\cite{HK00}] \label{BGP} 
Let $S$ be a spectrum whose mod $2$ cohomology is an injective unstable $\Ac$-module.  
Then $S$ satisfies the Brown-Gitler property if and only if either of the following condition holds:

\begin{itemize}
\item 
the evaluation 
map induces a monomorphism $H^* S  \ra H^*\Sigma^\infty \Omega^\infty S$ in cohomology,
\item  the spectrum is space-like,
\item there exists a space-like spectrum $Z$ 
and a map $ f : S 
\ra Z$ which induces  an epimorphism in cohomology. 
\end{itemize}
\end{theorem}

It is clear that a spectrum satisfying the Brown-Gitler property is uniquely defined by its cohomology. 

Recall that an unstable module $M$ is reduced, if either it does not contain a non-trivial suspension, or if the map from $M$ into itself, $x \in M^n$, $\Sq_0 \colon x \mapsto \Sq^{n}x$ is injective.

Below are the main examples of interest for us of such spectra.

\subsection{The Mitchell-Priddy  spectra}

The standard $\gl_n(\F)$ action  on $V_n$ induces one on $BV_n$, on the direct sum of all non-trivial line bundles over $BV_n$, and on the direct sum of $k$ copies of this vector bundle. Thus it acts also
on their Thom spaces
$BV_n^{\reg_n^{\oplus k}}$. 
As a consequence the Steinberg idempotent  $e_n$  \cite{Ste56} induces a stable map on 
$BV_n^{\reg_n^{\oplus k}}$, this map is defined after one suspension. 
Consider the telescope spectrum of this map  \cite{Tak99}, this is a stable summand of  $BV_n^{\reg_n^{\oplus k}}$, denoted by  $e_n\cdot BV_n^{\reg_n^{\oplus k}}$. Following Takayasu  we will shorten 
$e_n\cdot BV_n^{\reg_n^{\oplus k}}$ to  $\M(n)_k$.

In particular  $\M(n)_0=\M(n)$, $\M(n)_1=\L(n)$ are the spectra of Michell and Priddy and by definition  $\M(n)_2=\L'(n)$.  The spectra 
$\M(n)$ and $\L(n)$ satisfy the Brown-Gitler property.

The cohomology  of these spectra follows from the Thom isomorphism theorem:
$$H^*\M(n)_k=\omega_n^ke_n\F[x_1,\cdots,x_n]=\omega_n^{k-1}L_n$$
where  $\omega_n$ is the Euler class of the reduced regular representation $\reg_n$, or the top Dickson invariant of $\gl_n(\F)$: the product of all non-trivial linear forms in $\F[x_1,\cdots,x_n]$. 

\begin{theorem}[Takayasu \cite{Tak99}] For $k\ge 0$, there exist a cofibration sequence
$$\Sigma^k\M(n-1)_{2k+1}\xrightarrow{i_{n,k}} \M(n)_k\xrightarrow{j_{n,k}} \M(n)_{k+1}.$$
This induces the following short exact sequence in cohomology:  
$$0\rightarrow \omega_n^{k} L_n \rightarrow \omega_n^{k-1} L_n \rightarrow \Sigma^k\omega_{n-1}^{2k} L_{n-1}\rightarrow 0.$$
\end{theorem}
Splicing the cofibration sequences together yields:
{\small
\begin{equation*}\label{T}
\Sigma^{2^n-1}\M(0)_{2^n}\rightarrow \cdots 
\rightarrow \Sigma^{2^k-1} \M(n-k)_{2^k}\rightarrow \Sigma^{2^{k-1}-1} \M(n-k+1)_{2^{k-1}}\rightarrow
\cdots \rightarrow \Sigma \M(n-1)_2\rightarrow \M(n)_1\rightarrow \M(n)_2.
\end{equation*}
}

\subsection{The Brown-Gitler  spectra}

Let $J(\ell)$ the $\ell$-th Brown-Gitler module (cf. \cite[Chapter 2]{Sch94}).  It is characterised  by the fact that there is a natural equivalence 
$$
{\rm Hom}_\U(M,J(\ell)) \cong M^{\ell*}.
$$
The unstable module $J(\ell) $ is the injective hull of $\Sigma^\ell  \F$ in the category $\U$.
It is known that there  exist spectra $\T (\ell)$   
with   mod
  $2$ singular  cohomology  $J(\ell)$ 
\cite{GLM92}. This was originally proved by E. Brown et S. Gitler. If $\ell$ is odd $\T(\ell)$  is a suspension 
spectrum. More precisely we have

\begin{theorem} Let $n>2$. There exists a $1$-connected pointed space $\TT_1(n)$  (that we are going to call the $n$-th Brown-Gitler space),  unique up to homotopy,   such that the following holds: 
\begin{itemize}
\item There is an isomorphism of unstable modules $H^*\TT_1(n) \cong \Sigma J(n)$.
\item The induced map in cohomology by the evaluation   $\Sigma \Omega \TT_1(n) \ra \TT_1(n)$  is a monomorphism.
\item $\TT_1(n)$ is $2$-complete.
\end{itemize}

\noindent This space has moreover the following properties:
 \begin{itemize}
\item It is retract of the suspension of a pointed space.
\item For all pointed spaces $Y$ the natural map
$$
[\TT_1 (n), \Sigma Y]_* \ra  {\rm Hom}_\U(\tilde H^*(\Sigma Y), \Sigma J(n))
$$ is surjective.
\end{itemize}
\end{theorem}

Thus  the suspension spectrum of $\TT_1(n)$ is $\Sigma \TT(n)$, and  if $n$ is odd $\TT(n)$ is the suspension spectrum of the space  $\TT_1(n-1)$.

These were the first known spectra to have the two first conditions of Theorem  \ref{BGP}.

\subsection{The Lannes-Zarati spectra \cite{LZ86}}

\begin{theorem}
The unstable module $L_n \otimes J(k)$ is an injective object  in $\U$.
\end{theorem}

This is just a particular case of the Lannes-Zarati theorem. This module is the reduced mod $2$ cohomology of $\L(n) \wedge \TT(k)$. If $k$
is odd this is the suspension spectrum  of  $\L(n) \wedge \TT_1(k-1)$.  This is because $\TT_1(k-1)$ has a co-H-space structure.

\section{Proof of theorem \ref{NS}}
Let us recall some notations. 
The morphisms $$f_{n-k,n} \colon L_{k+1} \otimes J(2^{n-k-1}-1) \ra L_{k} \otimes J(2^{n-k}-1)$$ of the complex of Theorem
\ref{NST} are  defined in \cite{NST09b,NST09a}. Let us denote by $C_k$ the image
of the map $f_{n-k,n}$, which is also the kernel of $f_{n-k+1,n}$. 
If the spectra $\C(k)$ exist one has $H^*\C(k)\cong C_k $.

We are now going to proceed by descending induction. 
The first step of the proof is done by observing that the cofibration 
$$
\C(n-1) \ra \L(n)  \ra \C(n)
$$
is nothing else than Takayasu's cofibration 
$$
\Sigma \M(n-1)_3 \ra \M(n)_1\rightarrow \M(n)_2.
$$
However the right map is just the zero section, and the left hand spectrum could be defined as above.

We will assume that we have got a  spectrum  $\C(k-1)$, $2\le k\le n$, which is  the fibre of a map from $\L(k)\wedge \T(2^{n-k}-1)$ 
to  $\C(k)$  and  short exact sequences in cohomology $$H^*(\C(k)) \ra H^*(\L(k)\wedge \T(2^{n-k}-1))  \ra H^*(\C(k-1)).$$
 
We will make use of the following:
 
 \begin{proposition} [\cite{HK00}]  Let $X$ be a 
 $0$-connected spectrum such that $H_*X$ is of finite type and let $Z$ 
 be a space. Suppose there exists a map $f \colon X \ra \Sigma^\infty Z$ which is an epimorphism in cohomology. Then the evaluation map $ \Sigma^\infty  \Omega^\infty X \ra X$ induces a monomorphism in cohomology. 
\end{proposition}

One can apply this proposition to the spectrum $\C(k-1)$ and the suspension spectrum $\Sigma^\infty Z$ 
of which $\L(k) \wedge \T(2^{n-k}-1)$ is a wedge summand.  Thus the  evaluation map 
$\Sigma^\infty \Omega^\infty \C(k-1) \ra  \C(k-1)$ induces a monomorphism in cohomology:
$$
C_{k-1}  \cong  H^*\C(k-1) \ra H^*\Omega^\infty \C(k-1).
$$

As a consequence, as $L_k \otimes J(2^{n-k}-1)$ is injective in $\U$, the canonical map 
$C_{k-1} \ra L_{k-1} \otimes J(2^{n-k+1}-1)$ extends to a map of unstable modules 
$$H^*\Omega^\infty \C(k-1)\ra L_{k-1} \otimes J(2^{n-k+1}-1).$$
It follows 
from \ref{BGP} that this map is induced by a map of spectra 
$$\gamma \colon \L(k-1) \wedge \T(2^{n-k+1}-1) \ra \Sigma^\infty \Omega^\infty \C(k-1).$$ 
Composing $\gamma$ with the evaluation $\epsilon\colon\Sigma^\infty \Omega^\infty \C(k-1)
\ra \C(k-1)$  one gets a map 
$$\epsilon\gamma\colon \L(k-1)\wedge \T(2^{n-k+1}-1) \rightarrow \C(k-1)$$
inducing in cohomology the canonical map $C_{k-1} \ra  L_{k-1} \otimes J(2^{n-k+1}-1)$. 
One takes the fibre $\C(k-2)$ of the map $\epsilon\gamma$ 
to obtain a cofibration of spectra which induces a short exact sequence in cohomology.

The final step of the induction deserves care. One needs to show 
that the spectrum $\C(0)$ is indeed equivalent to $\T(2^n-1)$, {\it a priori} it has the right cohomology. By construction there is 
a map $\C(0) \ra \L(1)\wedge \T(2^{n-1}-1)$ which induces an epimorphism in cohomology. As $\L(1)\wedge \T(2^{n-1}-1)$ is space-like, it follows from 
\ref{BGP} that $\C(0)$ satisfies the Brown-Gitler property, thus 
$\C(0)\simeq\T(2^n-1)$. 
\begin{remark}
\begin{enumerate}
\item One has a fibration of spaces:
$$
\Omega^\infty \C(k-1) \ra \Omega^\infty \L(k-1) \wedge \T(2^{n-k+1}-1) \ra \Omega^\infty \C(k).
$$
This is even a fibration of H-spaces which allows us to compute $H^*\Omega^\infty \C(k-1)$. In fact the Eilenberg-Moore spectral sequence shows that:
$$
H^*\Omega^\infty \C(k-1) \cong H^* (\Omega^\infty \L(k-1) \wedge \T(2^{n-k+1}-1)) \otimes _{H^*\Omega^\infty \C(k)} \F.
$$
However we don't need this computation in the above proof, thus we don't describe (in terms of Dyer-Lashof operations) the result further.
\item The final step of the proof suggests it is possible to describe the space $\T_1(2^n-2)$ as a sub-complex of $(B\Z/2)^{\wedge n}$. We hope to be able to say more about that later.
\end{enumerate}
\end{remark}


\providecommand{\bysame}{\leavevmode\hbox to3em{\hrulefill}\thinspace}
\providecommand{\MR}{\relax\ifhmode\unskip\space\fi MR }
\providecommand{\MRhref}[2]{%
  \href{http://www.ams.org/mathscinet-getitem?mr=#1}{#2}
}
\providecommand{\href}[2]{#2}

\noindent LAGA, Universit\'e de Paris XIII, 93430 Villetaneuse,
France

\end{document}